\documentclass{article}

\usepackage{mathtools}
\usepackage{amsmath}
\usepackage{hyperref}
\usepackage{amsfonts}
\usepackage{xcolor}

\newcommand{\qed}{\mbox{}\nolinebreak\hfill\rule{2mm}{2mm}\par\medbreak}
\newcommand{\Proof}{\noindent{\bf Proof. }}

\usepackage[left=2.8cm,top=2cm,bottom=2cm,right=2.8cm]{geometry}

\usepackage{scrtime}
\usepackage[multiple]{footmisc}
\makeatletter
\renewcommand*\@fnsymbol[1]{\the#1}  
\makeatother

\def\RR{{ I\!\!R}}

 \newtheorem{proposition}{Proposition}
 \newtheorem{lemma}{Lemma}
 \newtheorem{theorem}{Theorem}
 \newtheorem{corollary}{Corollary}[theorem]
 
 \newtheorem{remark}{Remark}

\title{Consensus, polarization, and optimization of the mean value in a nonlinear model of opinion dynamics}

\author{David N. Reynolds\footnote{Partially supported by Modeling Nature (MNAT) Research Unit (project QUAL21-011)}, Pedro J.
Torres\footnote{Partially supported by project PN-Proyecto Plan Nacional MCIU-22-PID2021-128418NA-I00}\\
[4mm] {\small Departamento de Matem\'atica Aplicada \& }\\ {\small Modeling Nature (MNat) Research Unit.,}
\\ {\small Universidad de Granada, E-18071 Granada, Spain.}
}

\begin{document}
  \maketitle

\noindent {\bf Abstract}: This paper investigates some aspects of a recently proposed nonlinear mathematical model of opinion dynamics. The main objective is to identify the network structures that maximize the average equilibrium opinion (HMO). We prove that consensus is not generally attainable for populations with heterogeneous convictions, and that the highest mean does not necessarily correspond to consensus. Our analysis includes a necessary and sufficient condition for achieving the HMO, description of an algorithm for constructing optimal connectivity matrices, and strategies for pruning agents when heterogeneity obstructs mean optimization.\\


\noindent {\bf AMS2000 Subject
Classification: 91D30}.

\section{Introduction}\label{Seccion introduccion}

Public opinion is a driving force of paramount importance in modern societies, acting as a powerful feedback loop between citizens and decision-makers. The dynamic nature of opinion formation arises from continuous interactions among individuals, each influenced by the prevailing views of their peers as well as their own intrinsic convictions. Capturing the evolution of these opinions, especially in heterogeneous communities, remains a challenge of both theoretical interest and practical significance. In this regard, mathematical modeling of opinion dynamics has become an essential tool to understand how consensus, polarization, and collective behaviors arise within a group of individuals. 

Although the consideration of this topic from the point of view of mathematical modelling is relatively recent, many approaches have been proposed. Surveys such as \cite{H,P,PT1,PT2,SLST} offer valuable insights but remain necessarily incomplete due to the rapidly expanding body of related literature. In this paper, we concentrate on the recent nonlinear model,
\begin{equation}\label{model}
    x_{i}'(t)=\sum_{j=1}^{N}a_{ij}(x_{j}(t)-x_i(t))+ \sigma_{i}\left(u_i-x_{i}(t)\right)x_i(t),\quad i=1,\ldots,N.
\end{equation}
 This system is proposed in \cite{LRS,ReTu} to model the opinion dynamics of a community of $N$ agents, where $x_i(t)$ denotes the opinion value at time $t$, $u_i$ represents a constant personal conviction value, and each $\sigma_i$ models a fixed degree of each agent's stubbornness. The connectivity matrix $A=(a_{ij})$ is a symmetric matrix with nonnegative entries. Then, the first component rules the influence among agents, driving them toward consensus. On the other hand, the second component modulates the extent to which individuals adhere to their personal beliefs. Variants of this model have been studied in a different context, as a model for metapopulations in Mathematical Ecology \cite{RB,H,RuTo}, but its consideration as a model for opinion dynamics is relatively new, though related to classical linear models seen in \cite{BPDS,FJ}, as well as bounded confidence models \cite{Heg}, and flocking models \cite{CS2007a}.  

Denote by $u=(u_1,\ldots,u_N) \in \mathbb{R}^N_+$ the vector of \textit{conviction} values, $\sigma=(\sigma_1,\ldots,\sigma_N)\in \mathbb{R}^N_+$, the vector of \textit{stubbornness} parameters. The main result of \cite{ReTu} states that for any triplet $(A,u,\sigma)$, with $A$ symmetric, the system \eqref{model} has a unique equilibrium $x^*=(x_1^*,\ldots,x_N^*) \in \mathbb{R}^N_+$, which is a global attractor, and in general is representative of a \textit{compromise} state where for all $j$, $\min_j u_j<x_j<\max_j u_j$ and in general $x_j\neq x_k$ so that final opinions are closer together than the original conviction values, but not at consensus. The main problem we intend to address in this paper is the following question: for fixed values $u,\sigma$, {\bf which network connectivity matrix $A$ maximizes the average opinion at equilibrium?}. In the rest of the paper, the highest mean opinion state, if it exists, will be abbreviated by {\bf HMO}.

This question is relevant in many practical examples related with social networks, brand positioning strategies, political perceptions, etc. In the literature, there is a considerable number of papers devoted to strategies to reach a consensus state (a state in which all agents share the same opinion value), especially after the pioneering work \cite{DG} (see the review \cite{H}). In many practical situations related to group decision-making (like a jury or a committee), reaching a consensus is the desirable outcome, but as we will see later, such a state does not always coincide with the {\bf HMO} state.

Let ${\cal A}$ represent the set of symmetric matrices with non-negative components and null diagonal. Our objective, in mathematical terms, is to determine whether the function ${\cal M}: {\cal A}\to\RR^+$ given by
$$
{\cal M}(A)=\frac{1}{N}\sum_{i=1}^{N}x_i^*,
$$
reaches an absolute maximum, locate this maximum, (if possible) describe the set of connectivity matrices on which it is achieved, and (when not possible) prescribe a means for attaining the best possible reduced \textbf{HMO}.

To address the main question of identification of {\bf HMO} and its relation with consensus, we first analyze in Section \ref{sect2} the notion of consensus in our particular model, and how it can be asymptotically approached for highly connected systems. Section \ref{sect3} provides an upper bound for the {\bf HMO}. Then, in order to gain insight about the main question, we perform in Section \ref{sect4} a detailed analysis of the case of two agents. This particular case is the basis for a general study in Section \ref{sect5}, where a necessary and sufficient condition for the achievement of {\bf HMO} is given, including precise information of how to construct the optimal connectivity matrix. Some examples are included to illustrate the strategy. In Section \ref{sect6} we present the general case of uniform stubbornness $\sigma_i\equiv \sigma$, proving that the \textbf{HMO} is always attainable via the algorithm provided in Section \ref{sect5}. Section \ref{sect7} provides a notion of \textit{pruning} which is used to remove agents from the connectivity matrix in order to find the \textbf{HMO} when the algorithm in Section \ref{sect5} fails due to heterogeneity of stubbornness parameters. We finish the paper with a brief section of conclusions.

\section{Consensus for high connectivities}\label{sect2}

For fixed values of $u,\sigma$ and a given connectivity matrix $A\in {\cal A}$, the variance of the equilibrium $x^*$ is defined as
$$
{\cal V}(A)=\frac{1}{N}\sum_{i=1}^{N}(x_i^*-{\cal M}(A))^2.
$$
A perfect consensus is when all the coordinates of the equilibrium $x^*$ are equal. Then, the variance is a measure of how far a given system is from perfect consensus. Perfect consensus is achieved when ${\cal V}(A)=0$. The first result shows that it is not possible to achieve perfect consensus for a heterogeneous population.

\begin{proposition}\label{prop1}
 ${\cal V}(A)=0$ if and only if $u_i^*=u_j^*$ for every $i,j=1,\ldots,N$.
\end{proposition}

\Proof 
The non-trivial equilibrium satisfies the system of algebraic equations
\begin{equation}\label{syst}
    \sum_{j=1}^{N}a_{ij}(x_{j}^*-x_i^*)+ \sigma_{i}\left(u_i-x_{i}^*\right)x_i^*=0,\quad i=1,\ldots,N.
\end{equation}But ${\cal V}(A)=0$ means that $x_i^*=x_j^*$ for any $i,j=1,\ldots,N$. Then, the first term is zero and $x_i^*=u_i$ for any $i=1,\ldots,N$. The reverse implication is trivial
\qed

\begin{theorem}
    Assume that $A\in{\cal A}$ is irreducible. Then,
    $$
    \lim_{\lambda\to +\infty} {\cal V}(\lambda A)=0
    $$
    and
    $$
    \lim_{\lambda\to +\infty} {\cal M}(\lambda A)=\frac{\sum_{i=1}^{N}\sigma_{i}u_i}{\sum_{i=1}^{N}\sigma_{i}}.
    $$
\end{theorem}

\Proof The system verified by the equilibrium is 
\begin{equation}\label{Syst}
    \lambda\sum_{j=1}^{N}a_{ij}(x_{j}-x_i)+ \sigma_{i}\left(u_i-x_{i}\right)x_i=0,\quad i=1,\ldots,N,
\end{equation}
where we have dropped the $*$ for convenience. Adding the equations, we find 
\begin{equation}\label{rest}
  \sum_{i=1}^{N}\sigma_{i}\left(u_i-x_{i}\right)x_i=0.
\end{equation}
It is easy to see that this equation defines a compact set ${\cal K}$ of $\RR^N$, not depending on $\lambda$. Then, writing system \eqref{Syst} as
\begin{equation}\label{syst2}
\sum_{j=1}^{N}a_{ij}(x_{j}-x_i)+\frac{1}{\lambda} \sigma_{i}\left(u_i-x_{i}\right)x_i=0,\quad i=1,\ldots,N,
\end{equation}
and passing to the limit, the equilibrium when $\lambda\to+\infty$ leads to the linear homogeneous system holding
$$
\sum_{j=1}^{N}a_{ij}(x_{j}-x_i)=0,\quad i=1,\ldots,N.
$$
Call the corresponding matrix of coefficients $M$. Then, since the sum of rows (and columns) of $M$ is null, $M$ has $0$ as an eigenvalue with associated eigenvector $v=(1,\ldots,1)^T$. It remains to show that such eigenvalue is simple. To this aim, we observe that $M=A-D$, where $D$ is a diagonal matrix with positive elements in its diagonal (because of the irreducibility of A). Then, $v\in Ker(A)$ if and only if $Av=Dv$, i.e., $D^{-1}A v=v$. But the matrix $D^{-1}A$ is non-negative and irreducible, then by the Perron-Frobenius Theorem, $v$ is a simple eigenvector.

Based on the argument above, at the limit $\lambda\to+\infty$ the coordinates of the equilibrium have the same value, which of course is the mean value of them. This mean value, call it $m$, can be computed explicitly. Passing to the limit in the restriction \eqref{rest}, we have
$$
 \sum_{i=1}^{N}\sigma_{i}\left(u_i-m\right)m=0, 
$$
so that
$$
m=\frac{\sum_{i=1}^{N}\sigma_{i}u_i}{\sum_{i=1}^{N}\sigma_{i}}.
$$
This concludes the proof.
\qed
\section{Upper bound for the {\bf HMO}}\label{sect3}

\begin{theorem}\label{th1}
For any $A\in{\cal A}$, the following bound holds
$$
{\cal M}(A)\leq {\cal M}^*:=\frac{1}{2N}\left(\sqrt{\sum_{i=1}^{N}\sigma_i^{-1}}\sqrt{\sum_{i=1}^{N}\sigma_i u_i^2}+\sum_{i=1}^{N}u_i\right).
$$
\end{theorem}

\Proof
At equilibrium, we deal with the system of algebraic equations
\begin{equation}\label{syst}
    \sum_{j=1}^{N}a_{ij}(x_{j}-x_i)+ \sigma_{i}\left(u_i-x_{i}\right)x_i=0,\quad i=1,\ldots,N,
\end{equation}
where we have dropped the $*$ for convenience. Adding the equations, we find the restriction \eqref{rest}, which defines a compact set ${\cal K}$ of $\RR^N$. Thus, the function $F:(x_1,\ldots,x_N)\to \frac{1}{N}\sum_{i=1}^{N}x_i$ when restricted to ${\cal K}$ has a global maximum. To compute it, we use the standard method of Lagrange multipliers. Define
$$
L(x_1,\ldots,x_N,\lambda)=\frac{1}{N}\sum_{i=1}^{N}x_i+\lambda \sum_{i=1}^{N}\sigma_{i}\left(u_i-x_{i}\right)x_i.
$$
The critical points are easily obtained by solving
$$
\frac{\partial L}{\partial x_i}=\frac{1}{N}+\lambda\sigma_{i}\left(u_i-2x_{i}\right)=0,
$$
giving
\begin{equation}\label{equ}
x_i=\frac{1}{2N\lambda\sigma_i}+\frac{u_i}{2}.
\end{equation}
Now, $\lambda$ is obtained by substituting the latter values of $x_i$ in the restriction \eqref{rest}. Some easy computations provide
\begin{equation}\label{lambda}
\lambda=\frac{1}{N}\sqrt{\frac{\sum_{i=1}^{N}\sigma_i^{-1}}{\sum_{i=1}^{N}\sigma_i u_i^2}}
\end{equation}
Therefore,
$$
{\cal M}^*=\frac{1}{N}\sum_{i=1}^{N}\left(\frac{1}{2N\lambda\sigma_i}+\frac{u_i}{2}\right)=\frac{1}{2N}\left(\sqrt{\sum_{i=1}^{N}\sigma_i^{-1}}\sqrt{\sum_{i=1}^{N}\sigma_i u_i^2}+\sum_{i=1}^{N}u_i\right).
$$
\qed

\section{Analysis of the case $N=2$}\label{sect4}
For the case of two agents, the set of connectivity matrices is the one-parameter family 
$$
A=\left(\begin{array}{cc}
 0  & a \\
a & 0
\end{array}\right)
$$
with $a\geq 0$. Hence, by convenience, we will identify ${\cal M}(A)\equiv {\cal M}(a)$.

\begin{lemma}
Assume that $u_1\leq u_2$ (this can be done without loss of generality). Then,
\begin{equation}\label{ineq}
u_1\leq x_1^*\leq x_2^*\leq u_2.
\end{equation}

\end{lemma}

\Proof 

At equilibrium, the system is

\begin{equation}\label{syst-N=2}
\left\{\begin{array}{ccc}
 a(x_2^*-x_1^*)  &=& \sigma_1(x_1^*-u_1)x_1^*, \\
a(x_1^*-x_2^*) &=& \sigma_2(x_2^*-u_2)x_2^*.
\end{array}\right.
\end{equation}
Arguing by contrapositive, assume that $x_1^*>x_2^*$, then from the first equation we would have $x_1^*<u_1$, and from the second equation, $x_2^*>u_2$. Hence, $u_1>u_2$. In consequence, it is proved that $x_1^*\leq x_2^*$. Now, the whole sequence of inequalities is easily proved from \eqref{syst-N=2}.
\qed

\begin{remark}
    A direct consequence of this lemma is that if $u_1=u_2$, then the equilibrium is a consensus $x_1^*=x_2^*=u_1$, which is consistent with Proposition \ref{prop1}. 
\end{remark}

Our objective is to find the maximum of ${\cal M}(a)$, if it exists. To this aim, the parameter
$$
\mu=\frac{\sigma_1}{\sigma_2}
$$
will play a fundamental role.

\begin{theorem}\label{trichotomy} Assume that $u_1 < u_2$. The subsequent statements are valid
\begin{itemize}
\item[(i)] (Polarization) If $\mu\geq \frac{u_2^2}{u_1^2}$, then
 ${\cal M}(a)< {\cal M}(0)=\frac{u_1+u_2}{2},$ for any $a>0$.
 
\item [(ii)] (Consensus) There exists a value $\mu_*>0$ such that if $\mu\leq \mu_*$, then  for any $a>0$,
$$
{\cal M}(a)< {\cal M}(+\infty):=\lim_{a\to +\infty}{\cal M}(a)=\frac{\sigma_1 u_1+\sigma_2 u_2}{\sigma_1+\sigma_2}
$$
\item[(iii)] (Compromise)  If $\mu_*<\mu< \frac{u_2^2}{u_1^2}$, there exists an explicitly computable $a^*>0$ such that 
$${\cal M}(a^*)=\max_{0<a<+\infty}{\cal M}(a).$$

\end{itemize}

\end{theorem}

\Proof
Adding the equations of system \eqref{syst-N=2}, the equilibrium $(x_1^*,x_2^*)$ belongs to the ellipse $\cal E$ of equation
$$
\sigma_1(x_1-u_1)x_1+\sigma_2(x_2-u_2)x_2=0.
$$
Moreover, it must verify \eqref{ineq}. Therefore, we are compelled to find the maximum of the function
$$
F(x_1,x_2)=\frac{x_1+x_2}{2}
$$
resctricted to the piece of ellipse $\cal E$ that verifies \eqref{ineq}. The maximum can be attained in the interior or in one of the extremes of the piece of ellipse. The extremes correspond to $(u_1,u_2)$ and the consensus $(\tilde x_1,\tilde x_1)$, with 
$$
\sigma_1(\tilde x_1-u_1){\tilde x}_1+\sigma_2({\tilde x}_1-u_2){\tilde x}_1=0,
$$
giving
$$
{\tilde x}_1=\frac{\mu u_1+u_2}{\mu+1}.
$$
The maximum of $F$ restricted to $\cal E$ was computed in the proof of Theorem \ref{th1}. Considering \eqref{equ}-\eqref{lambda}, it is attained at the point
\begin{equation}\label{x1-x2}
x_1=\frac{u_1}{2}+\frac{1}{2}\sqrt{\frac{\mu u_1^2+u_2^2}{\mu(1+\mu)}}, \qquad x_2=\frac{u_2}{2}+\frac{1}{2}\sqrt{\frac{\mu(\mu u_1^2+u_2^2)}{1+\mu}}.
\end{equation}
Now, we have to check if $(x_1,x_2)$ verifies \eqref{ineq} or not. After some trivial algebra, it turns out that $x_1\leq u_1$ if and only if $\mu\geq \frac{u_2^2}{u_1^2}$. In this case, the maximum is attained at the extreme $(u_1,u_2)$  of the piece of the ellipse and $(i)$ holds. 

Conversely, $x_1<x_2$ means that
\begin{equation}\label{ineq-2}
    u_1+\sqrt{\frac{\mu u_1^2+u_2^2}{\mu(1+\mu)}}<u_2+\sqrt{\frac{\mu(\mu u_1^2+u_2^2)}{1+\mu}}.
\end{equation}
The left-hand side of this inequality tends to $+\infty$ when $\mu\to 0^+$, whereas the right hand side tends to $u_2$. This means that there exists a value $\mu_*>0$ such that $x_1>x_2$ if $\mu<\mu_*$. In this case, the maximum is in the consensus extreme $(\tilde x_1,\tilde x_1)$, which is approached asymptotically as $a\to +\infty$. Therefore, the point $(ii)$ is proved.

Finally, to prove $(iii)$, note that if $\mu_*<\mu< \frac{u_2^2}{u_1^2}$, then the point $(x_1,x_2)$ defined by \eqref{x1-x2} satisfies \eqref{ineq}. Then, we can define
$$
a^*=\frac{\sigma_1(x_1-u_1)x_1}{x_2-x_1}>0.
$$
By construction, ${\cal M}(a^*)\geq {\cal M}(a)$ for all $a>0$, and the proof is finished.

\qed

\begin{remark}
    The threshold value $\mu_*$ can be computed numerically in any concrete case, but it is easy to give an explicit estimate: note that the inequality \eqref{ineq-2} holds for $\mu=1$, which means that $ \mu_*\leq 1$.
\end{remark}

\section{Identification of the optimal connectivity matrix}\label{sect5}

The proof of Theorem \ref{th1} provides a piece of important information: if the maximum is attained, the equilibrium must be located at the coordinates given by \eqref{equ}, with $\lambda$ given by \eqref{lambda}. In this way, the `ideal' equilibrium is identified explicitly in terms of the known parameters $u,\sigma$.

Let us call
\begin{equation}\label{equ-virt}
x_i^*=\frac{1}{2N\lambda\sigma_i}+\frac{u_i}{2},
\end{equation}
where $\lambda$ defined by \eqref{lambda}. Then, we deal with a kind of inverse problem, that is to look for a suitable matrix $A\in{\cal A}$ such that 
\begin{equation}\label{syst-2}
    \sum_{j=1}^{N}a_{ij}(x_{j}^*-x_i^*)+ \sigma_{i}\left(u_i-x_{i}^*\right)x_i^*=0,\quad i=1,\ldots,N,
\end{equation}
To simplify the writing, we define
$$
d_{ij}=x_{j}^*-x_{i}^*, \quad g_i=\sigma_{i}\left(x_{i}^*-u_i\right)x_i^*.
$$
Then, \eqref{syst-2} is a linear system 
\begin{equation}\label{syst-3}
    \sum_{j=1}^{N}d_{ij}a_{ij}= g_i,\quad i=1,\ldots,N,
\end{equation}
with unknowns $a_{ij}$. Observe that the sum of the equations vanishes. Since $a_{ij}=a_{ji}$ and $a_{ii}=0$, this system has $N(N-1)/2$ unknowns and $N$ equations, so it will typically be indeterminate compatible. The problem is to determine if there are solutions with non-negative values. The reference \cite{Dines} provides a constructive algorithm to solve this problem in general. We develop here a tailored strategy.

For simplicity, let us consider the case where all the coordinates given by \eqref{equ-virt} are different. Then, without loss of generality, we can reindex the agents to order such scalars as
$$
x_1*<x_2^*<\ldots<x_N^*.
$$
This means that $d_{ij}>0$ for all $i<j$.

\begin{theorem}\label{algorithm}
    The condition
    \begin{equation}\label{suf}
\sum_{i=1}^{k} g_i>0,\quad k=1,\ldots,N-1,
    \end{equation}    
is necessary and sufficient for the existence of a non-trivial non-negative solution of \eqref{syst-3}.  
\end{theorem}

\Proof 
Assume that \eqref{syst-3} has a positive solution $a_{ij}$. Considering that $a_{ij}=a_{ji}$, $a_{ii}=0$, $d_{ij}=-d_{ji}$ and $d_{ij}>0$ for all $i<j$, condition \eqref{suf} follows from adding the $k$ first equations of the system for $k=1,\ldots,N-1$.

Assume now that \eqref{suf} holds. The first equation of the system is 
$$
d_{12} a_{12}+d_{13} a_{13}+\cdots+d_{1N} a_{1N}= g_1.
$$
Then, we fix
$$
a_{12}=\frac{g_1}{d_{12}}>0,\quad a_{1j}=0\quad j=3,\ldots,N.
$$
The second equation is 
$$
-d_{12} a_{12}+d_{23} a_{23}+\cdots+d_{2N} a_{2N}= g_2,
$$
which, considering the fixed value of $a_{12}$, can be written as
$$
d_{23} a_{23}+\cdots+d_{2N} a_{2N}= g_1+g_2.
$$
Then, as a solution we fix
$$
a_{23}=\frac{g_1+g_2}{d_{23}}>0,\quad A_{2j}=0\quad j=3,\ldots,N.
$$
By recurrence, we construct the whole connectivity matrix as
$$
a_{k,k+1}=\frac{1}{d_{k,k+1}}\sum_{i=1}^{k} g_i>0,\quad \quad a_{kj}=0,\quad k=1,\ldots,N-1,\,j=k+2,\ldots,N.
$$\qed

\subsection{N=2}

We analyze the simplest case of two agents. Since $d_{ij}=-d_{ji}$ and the connectivity matrix $A$ must be symmetric, system \eqref{syst-3} is just

\begin{equation}
\left\{\begin{array}{ccc}
d_{12} a_{12}  &=& g_1, \\
-d_{12} a_{12} &=& g_2.
\end{array}\right.
\end{equation}
Note that $g_1+g_2=0$. Then, 
$$
a_{12}=\frac{g_1}{d_{12}}.
$$

\noindent{\bf Examples. }
Let us provide three illustrative examples for the three possible outcomes of the $N=2$ case.

 \begin{enumerate}
\item (Polarization)   For $u=(1,2)$, $\sigma=(3,1)$, we get $x_1^*\sim 0.88$ and $x_2^*\sim 2.145$ $\implies$ $a_{1,2}\simeq -0.247277$ and that no positive connection can attain the \textbf{HMO} so $x_1^*=1,x_2^*=2$ and $\mathcal{M}^*=\frac{3}{2}$ is the highest attainable average opinion.
 
 \item (Consensus)  For $u=(1,5)$, $\sigma=(1,10)$, we get $x_1^*\sim 8.05$ and $x_2 \sim 3.25$ $\implies$ $a_{1,2}\simeq -11.8$ and that no positive connection can attain the \textbf{HMO} so the connection should be strengthened to get the consensus maximal average opinion at $\mathcal{M}^*=\frac{51}{11}$.
 
 \item (Compromise) For $u=(1,2)$, $\sigma=(1,1)$, we get $a_{1,2}=\frac{3}{4}$ as the optimal connectivity, and the opinion average is  ${\cal M^*}=\frac{3+\sqrt{10}}{4}\simeq 1.54057$.

\end{enumerate}
Note that the way in which the algorithm fails determines the next best way to achieve the \textbf{HMO}. In the polarization example, $g_1<0$ broke the algorithm, while for the consensus case $d_{12}<0$ is what cause the algorithm to fail. 
\subsection{N=3}

In this case, system \eqref{syst-3} is 
\begin{equation}
\left\{\begin{array}{ccc}
d_{12} a_{12}+d_{13} a_{13}&=& g_1, \\
-d_{12} a_{12}+d_{23} a_{23} &=& g_2, \\
-d_{13} a_{13}-d_{23} a_{23} &=& g_3.
\end{array}\right.
\label{syst-N=3}
 \end{equation}
To understand when this system admits solutions in the first orthant, we assume without loss of generality that $x_1^*\leq x_2^*\leq x_3^*$ (if it is not the case, we can always re-order the indices). Then, $d_{ij}\geq 0$ for every $i<j$. From the first and third equations, the condition $g_1>0>g_3$ is necessary. It turn out that it is also sufficient for the existence of a solution with positive entries.
 
\noindent{\bf Example. }
 For $u=(1,4,5)$, $\sigma=(1,2,1)$, the solution of \eqref{syst-N=3} is
 $$
\left\{\begin{array}{l}
a_{12}=13.5214+5.30249\mu, \\
a_{13}= 0.1893-0.841333\mu, \\
a_{23}=\mu.
\end{array}\right.
$$
For $\mu>0$ small enough, these coefficients are positive, and in particular, the value of $\mu$ that makes $a_{13}=0$ provides the connectivity matrix one would find using the algorithm from Theorem \ref{algorithm}. In this case, ${\cal M^*}\simeq 3.6736$.

\subsection{Optimal network for a large community with two subgroups}

In \cite{ReTu}, numerical simulations were done for a community of 150 agents with two subgroups with conviction and stubbornness vectors given by
$$
u_i=\begin{cases}
\kappa &  \mbox{ if }\quad 1\leq i\leq 50, \\
1 & \mbox{ if }\quad 51\leq i\leq 150,
\end{cases}
\qquad \mbox{ and }\qquad \sigma_i=\begin{cases}
\delta & \mbox{ if } \quad 1\leq i\leq 50, \\
1 & \mbox{ if }\quad 51\leq i\leq 150,
\end{cases}
$$
for different types of connectivity networks. The simulations therein highlighted the nonlinear effects of the stubbornness term $\sigma_i(u_i-x_i)x_i$. Our objective here is to identify the optimal network that may maximize the average opinion of the community. 

Following the notation of the beginning of the section, it is easy to compute the `ideal' equilibrium as
$$
x_i^*=\begin{cases}
\frac{\kappa}{2}+\frac{1}{300\delta} \sqrt{\frac{2\delta+\delta^2\kappa^2}{1+2\delta}}& \mbox{ if } \quad 1\leq i\leq 50, \\
 \frac{1}{2}+\frac{1}{300} \sqrt{\frac{2\delta+\delta^2\kappa^2}{1+2\delta}}& \mbox{ if }\quad 51\leq i\leq 150.
\end{cases}
$$
Hence, $d_{ij}=0$ if $1\leq i,j\leq 50$ or $51\leq i,j\leq 150$, while $d_{ij}=d$ (a fixed computable value) if $1\leq i\leq 50$ and $51\leq j\leq 150$. Furthermore,
$$
g_i=\begin{cases}
-2g & \mbox{ if }\quad 1\leq i\leq 50, \\
 g & \mbox{ if } \quad 51\leq i\leq 150.
\end{cases}
$$
where $g=\frac{\delta^2\kappa^2-1}{4+8\delta}$. Then, the system to be solved, corresponding to \eqref{syst-3}, is 
\begin{equation}\displaystyle
\begin{array}{ccc}
\displaystyle d\sum_{j=51}^{150}a_{ij}= -2g,\quad i=1,\ldots,50, \\
\displaystyle-d\sum_{j=1}^{50}a_{ij}= g,\quad i=51,\ldots,150.
\end{array}
\label{syst-N=150}
 \end{equation}
Here, we have used again that $d_{ij}=-d_{ji}$. By direct substitution, it is clear that a solution of this system is
$$
a_{ij}=\begin{cases}
\frac{g}{-50d} & \mbox{ if }\quad 1\leq i\leq 50 \mbox{ and } 51\leq j\leq 150 \mbox{ or } 1\leq j\leq 50 \mbox{ and } 51\leq i\leq 150 \\
 0 & \mbox{ elsewhere }.
\end{cases}
$$
Of course, since the system is compatible indeterminate, there are infinitely many other solutions, but this one is certainly the simplest structure. The interpretation is clear: to achieve the {\bf HMO}, the elements of the two different subgroups must be connected in a fixed way, while elements of the same subgroup do not need to be connected. We take the concrete values of parameters used in the numerical computations of \cite{ReTu} to compare results:
\begin{itemize}
\item $\kappa=100$, $\delta=2$: ${\cal M}^*=54.2697$ for the network
$$
a_{ij}=\begin{cases}
0.810502 & \mbox{ if }\quad 1\leq i\leq 50 \mbox{ and } 51\leq j\leq 150 \mbox{ or } 1\leq j\leq 50 \mbox{ and } 51\leq i\leq 150 \\
 0 & \mbox{ elsewhere }.
\end{cases}
$$

\item $\kappa=100$, $\delta=0.5$: ${\cal M}^*=40.5749$ for the network
$$
a_{ij}=\begin{cases}
0.125912 & \mbox{ if }\quad 1\leq i\leq 50 \mbox{ and } 51\leq j\leq 150 \mbox{ or } 1\leq j\leq 50 \mbox{ and } 51\leq i\leq 150 \\
 0 & \mbox{ elsewhere }.
\end{cases}
$$

\item $\kappa=10$, $\delta=2$: ${\cal M}^*=5.74537$ for the network
$$
a_{ij}=\begin{cases}
0.0889628 & \mbox{ if }\quad 1\leq i\leq 50 \mbox{ and } 51\leq j\leq 150 \mbox{ or } 1\leq j\leq 50 \mbox{ and } 51\leq i\leq 150 \\
 0 & \mbox{ elsewhere }.
\end{cases}
$$

\item $\kappa=10$, $\delta=0.5$: ${\cal M}^*=4.4037$ for the network
$$
a_{ij}=\begin{cases}
0.0132978 & \mbox{ if }\quad 1\leq i\leq 50 \mbox{ and } 51\leq j\leq 150 \mbox{ or } 1\leq j\leq 50 \mbox{ and } 51\leq i\leq 150 \\
 0 & \mbox{ elsewhere }.
\end{cases}
$$
In a natural way, higher stubbornness of the group with higher conviction value provides a higher average.

\section{A group of $N$ agents with equal stubbornness}\label{sect6}

We consider model \eqref{model} with  $\sigma_i\equiv\sigma$ for every $i=1\ldots N$. Then, 
$$
{\cal M}(A)= {\cal M}^*=\frac{1}{2N}\left(\sqrt{N}\sqrt{\sum_{i=1}^{N} u_i^2}+\sum_{i=1}^{N}u_i\right)
$$
 for a suitable choice of $A$, which we will show always exists.    
\end{itemize}

Let us rewrite the model with uniform stubbornness:
\begin{equation*}\label{uniformmodel}
    x_{i}'(t)=\sum_{j=1}^{N}a_{ij}(x_{j}(t)-x_i(t))+ \sigma\left(u_i-x_{i}(t)\right)x_i(t),\quad i=1,\ldots,N.
\end{equation*}
Then we see that fixed points must lie on the $N$-sphere defined by
\begin{align*}
    \sigma\sum_{i=1}^N(u_i-x_i)x_i=0,
\end{align*}
and the equations \eqref{equ} and \eqref{lambda} reduce to
\begin{equation}\label{equ-unif}
    x_i=\frac{1}{2N\lambda \sigma}+\frac{u_i}{2}
\end{equation}
\begin{equation}\label{lambda-unif}
    \lambda=\frac{1}{\sigma\sqrt{{N\sum_{i=1}^Nu_i^2}}}
\end{equation}

Therefore, the following lemma holds
\begin{lemma}
    If $u_1\leq \ldots \leq u_n$, then the optimal configuration satisfies
    \begin{align*}
        u_1\leq x_1\leq \ldots x_j\leq x_{j+1}\leq \ldots \leq x_N \leq u_N.
    \end{align*}
\end{lemma}
\Proof
    The inequality $x_j\leq x_{j+1}$ for each $j=1 \ldots N-1$ holds by \eqref{equ-unif} and the fact that $u_1\leq \ldots \leq u_n$. To see the confinement to the interval $[u_1,u_n]$ one must show that 
    \begin{align*}
        u_1\leq \frac{1}{N\lambda\sigma}\leq u_n.
    \end{align*}
    This holds from the fact that $u_1^2\leq \frac{1}{N}\sum_{j=1}^N u_j^2\leq u_N^2$ and \eqref{lambda-unif}.
\qed
Further, we can compute via the algorithm found in the proof of Theorem \ref{algorithm} an optimal connectivity matrix and further see that condition \eqref{suf} is guaranteed to hold.

\begin{theorem}\label{t:uniform}
    Suppose $\sigma_i\equiv \sigma>0$ and $0<u_i<u_{i+1}$ for each $i=1\ldots N-1$. Then the compromise state given by \eqref{equ-unif} with $\lambda$ satisfying \eqref{lambda-unif} gives an HMO state with the connectivity matrix
    \begin{align}
        a_{k,k+1}&=\frac{\sigma}{2}\left(\frac{\frac{k}{N}\sum_{i=1}^Nu_i^2-\sum_{i=1}^ku_i^2}{u_{k+1}-u_k}\right), \ \ \ k=1\ldots N-1, \label{HMOmat}\\
        a_{k,j}&=0, \ \ \ \text{if} \ \ j\neq k+1,
    \end{align}
    where $a_{k,k+1}>0$.
\end{theorem}
\Proof
Following the algorithm from Theorem \ref{algorithm}, we compute
\begin{align}
    d_{k,k+1}&=x_{k+1}-x_k=\frac{1}{2}(u_{k+1}-u_k), \nonumber\\
    g_k&=\sigma(x_k-u_k)x_k=\sigma\left(\frac{1}{4\lambda^2N^2\sigma^2}-\frac{u_k^2}{4}\right). \label{geq}
\end{align}
continuing with \eqref{geq} we can simplify with the equation \eqref{lambda-unif} to get
\begin{align*}
    g_k=\frac{\sigma}{4}\left(\frac{k}{N}\sum_{i=1}^Nu_i^2-\sum_{i=1}^ku_i^2\right).
\end{align*}

As by construction $a_{k,k+1}=\frac{1}{d_{k,k+1}}\sum_{i=1}^k g_i$ we have the formula \eqref{HMOmat} and we need only prove that $a_{k,k+1}>0$ for each $k=1 \ldots N-1$. As $0<u_k<u_{k+1}$ we need only prove the following
\begin{align*}
    \frac{k}{N}\sum_{i=1}^Nu_i^2-\sum_{i=1}^ku_i^2>0
\end{align*}
Since $\{u_k\}_{k=1}^N$ is a strictly increasing positive sequence this must hold due to the monotonicity of the partial arithmetic means. Indeed, letting $S_k=\sum_{i=1}^k u_i^2$ and $T_k=\sum_{i=k+1}^N u_i^2$ we have
\begin{align*}
    \frac{j}{N}\sum_{i=1}^Nu_i^2-\sum_{i=1}^ju_i^2=\frac{k-N}{N}S_k+\frac{k}{N}T_k.
\end{align*}
Further letting $m_1=\frac{S_k}{k}$ and $m_2=\frac{T_k}{N-k}$ we get,
\begin{align*}
    \frac{k}{N}\sum_{i=1}^Nu_i^2-\sum_{i=1}^ku_i^2=\frac{k(N-k)}{N}(m_2-m_1).
\end{align*}
Now as $m_1$ is the average of the first $k$ terms of $\{u_i^2\}_{i=1}^N$ and $m_2$ is the average of the last $N-k$ terms we have
\begin{align*}
    m_1\leq u_k^2<u_{k+1}^2\leq m_2,
\end{align*}
and hence $a_{k,k+1}>0$.

\qed 
Note that the above reinforces the fact that if $u_i=u_j$, then a link between agent $i$ and agent $j$ is unnecessary for finding the \textbf{HMO}, as was seen in the previous section. Further we see that under uniform stubbornness parameter $\sigma_i\equiv \sigma$, that the \textbf{HMO} is achieved at the compromise state on the interior of the $N$-sphere rather than the end points representing polarization and consensus states, respectively.

\section{Pruning the connectivity matrix}\label{sect7}

As one can see from part (i) of Theorem \ref{trichotomy}, if there is a discrepancy between stubbornness levels compared to conviction values strong enough, then the highest average opinion value is actually achieved by disconnecting the agents. We will refer to this process as pruning. For the $N=2$ case there is only one edge and so the trichotomy of Theorem \ref{trichotomy} states that the relationship between $\mathbf{\sigma}$ and $\mathbf{u}$ determines when the \textbf{HMO} is achieved, with a disconnected matrix (i) (polarizing), an infinitely strong matrix (ii) (consensus), or at a finite value (iii) (compromise). When $N>2$ we can obtain the \textbf{HMO} via combinations of these states.

\subsection{How to prune}
Reindex the variables so that $x_1^*< \ldots < x_N^*$ and recall the optimal values computed previously
\begin{align*}
    x_j^*=\frac{1}{2N\lambda\sigma_j}+\frac{u_j}{2}.
\end{align*}
If these optimal values are to be a fixed point of the model, then the following must hold
\begin{align}\label{min/max}
    0=\sum_{j=1}^Na_{1j}(x_j^*-x_1^*)+\sigma_1(u_1-x_1^*)x_1^*
\end{align}
and since $x_1^*$ is the minimal opinion value by the reindexing, we have $x_1^*\geq u_1$. Similarly we have $x_N^*\leq u_N$ must hold for the optimal values to be representative of an admissible fixed point. It is important to note that $u_1$ is not necessarily the minimal conviction value (and similarly for $u_N$ to be maximal).

Let us start with a first lemma.
\begin{lemma}\label{max/minrange}
    Let $0<x_1^*< \ldots <x_N^*$. Then the following must hold in order to be a fixed point of \eqref{model},
    \begin{align*}
        u_1\leq x_1^*< \ldots < x_N^*\leq  u_n,
    \end{align*}
    where again $u_1$ is not necessarily the minimal conviction value due to the reindexing (and similarly for $u_N$ to be maximal).
\end{lemma}
\Proof
Since equation \eqref{min/max} must hold at a fixed point, given that $0<x_1^*< \ldots <x_N^*$ we have $\sum_{j=1}^Na_{1j}(x_j^*-x_1^*)\geq0$ where $\sum_{j=1}^Na_{1j}(x_j^*-x_1^*)=0$ only occurs if $a_{1j}=0$ for all $j=1 \ldots N$.  Thus $x_1^*\geq u_1$. The same holds for $x_N^*\leq u_N$.
\qed

Note that the above does not show an ordering of the $x_j$ that is dependent on their corresponding $u_j$. However, it implies that the fixed point lies between the $u_1$ and $u_N$ that correspond to the minimal and maximal $x_j^*$, respectively.

The method for pruning comes from modifying the algorithm constructed in Section \ref{sect5}. Theorem \ref{algorithm} provided a necessary and sufficient condition for the maximal opinion value to be achieved with a nonnegative connectivity matrix. Clearly, when condition \eqref{suf} is met there is no need to prune. However, when there is heterogeneity of stubbornness parameters so that $\sigma_j \neq \sigma_k$, it is possible for the optimal opinion values to fall outside of the permissible range for a fixed point.

Indeed, as
\begin{align*}
    x_1^*=\frac{1}{2N\lambda\sigma_1}+\frac{u_1}{2}
\end{align*}
if $\sigma_1$ is large enough, then $x_1^*<u_1$ contradicting Lemma \ref{max/minrange}. In other words, during the algorithm for constructing the matrix one would get
\begin{align*}
    a_{12}=\frac{g_1}{d_{12}}<0
\end{align*}
since $g_1=\sigma_1(x_1^*-u_1)x_1^*<0$ and $d_{12}=x_2^*-x_1^*>0$ by construction. As negative connections are not allowed, a naive approach would be to prune $x_1$ from the system via setting $a_{1j}=0$ for all $j=1,...,N$. The system would then read
\begin{align*}
    x_i'&=\sum_{j=2}^N a_{ij}(x_j-x_i)+\sigma_i(u_i-x_i)x_i, \ \ i=2,...,N,\\
    x_1'&=\sigma_1(u_1-x_1)x_1.
\end{align*}
In this way we have $x_1 \to u_1$ exponentially fast and we have a new system of $N-1$ agents interacting so that the HMO is achieved at
\begin{align}
    x_1^*&=u_1,\label{psol1}\\
    x_j^*&=\frac{1}{2(N-1)\lambda_{p_1}\sigma_j}+\frac{u_j}{2},\label{psol2}\\
    \lambda_{p_1}&=\frac{1}{N-1}\sqrt{\frac{\sum_{i=2}^N\sigma_i^{-1}}{\sum_{i=2}^N\sigma_iu_i^2}}\label{psol3}
\end{align}
The removal of this agent from the connectivity matrix is known as \textit{pruning}.

Now if condition \eqref{suf} holds on the pruned network, then the algorithm yields a nonnegative connectivity matrix with unique fixed point at \eqref{psol1}-\eqref{psol3}. On the other hand, it is not clear, that pruning $x_1$ instead of another $x_j$ has the optimal effect on the final state's HMO.

Consider the following example,
\begin{itemize}
    \item $u=(1,2,3), \ \sigma=(1,1,0.1)$. Computing
    \begin{align*}
        x_1^*\sim 0.85, \ \ x_2^*\sim 1.35, \ \ x_3^*\sim 5,
    \end{align*}
    From Lemma \ref{max/minrange} we see that both $x_1^*<u_1$ and $x_3^*>u_3$ provide a negative conection in the connectivity matrix, and thus to achieve a system with a nonnegative solution we must prune. \begin{itemize}
        \item (Bottom pruning) If like above we prune $x_1$ we must recalculate the \textbf{HMO} state and we get
    \begin{align*}
        x_{1,p}^*=1, \ \ x_{2,p}^*\sim 1.33, \ \ x_{3,p}^*\sim 4.84,
    \end{align*}
    and since $x_{2,p}^*<u_2$ we must prune once more and the graph becomes disconnected giving a polarized state $x_p^*=(1,2,3)=u$ with $\mathcal{M}(A)=2$.
    \item (Top pruning) On the other hand if we prune $x_3$ from the network first we get
    \begin{align*}
        x_{1,p}^*\sim 1.29, \ \ x_{2,p}^*\sim 1.79, \ \ x_{3,p}^*=3,
    \end{align*}
    which gives
    \begin{align*}
        a_{12}=\frac{3}{4}, \ \ a_{13}=a_{23}=0,
    \end{align*}
    and $\mathcal{M}^*(A)\sim 2.026$.
    \end{itemize}
    Hence the top pruning led to a greater \textbf{HMO} achieved by pruning the top element from the network.
\end{itemize}

Heuristically it makes sense that agents with conviction value $u_i$ small, but stubbornness value $\sigma_i$ large detract from the overall average and should be pruned. Analogously the same can be said for $u_i$ large with $\sigma_i$ small. We will now formalize these heuristics for the case of one agent with varying stubbornness parameter and the other $N-1$ fixed.
\begin{lemma}\label{parameter}
    Let $x_i^*$ be the optimal state found earlier so that
    \begin{align*}
        x_i^*=\frac{1}{2N\lambda \sigma_i}+\frac{u_i}{2}.
    \end{align*}
    Then the following parametric behavior holds,
    \begin{align*}
        &\frac{dx_i^*}{d\sigma_i}<0, \ \ \ \lim_{\sigma_i \to \infty} x_i^*(\sigma_i)=\frac{u_i}{2}, \ \ \lim_{\sigma_i \to 0} x_i^*(\sigma_i)=\infty,\\
        &\lim_{\sigma_i \to \infty} \lambda(\sigma_i)=0, \ \  \lim_{\sigma_i \to 0} \lambda(\sigma_i)=\infty,\\
        &\lim_{\sigma_i\to \infty} x_j^*(\sigma_i)=\infty, \ \ \lim_{\sigma_i\to 0} x_j^*(\sigma_i)=\frac{u_j}{2}.
    \end{align*}
\end{lemma}
\Proof
As $u_i$ is fixed, we compute
\begin{align*}
    \frac{dx_i^*}{d\sigma_i}=\frac{-1}{2N\lambda^2\sigma_i^2}\frac{d(\lambda \sigma_i)}{d\sigma_i}.
\end{align*}
Computing first
\begin{align*}
    \frac{d\lambda}{d\sigma_i}=\frac{1}{2N^2\lambda }\left(\frac{\frac{1}{\sigma_i^2}-u_i^2N^2\lambda^2}{\sum_{j=1}^N\sigma_ju_j^2}\right).
\end{align*}
Thus
\begin{align*}
    \frac{d(\lambda \sigma_i)}{d\sigma_i}=\sigma_i\frac{d\lambda}{d\sigma_i}+\lambda=\frac{1}{2N^2\lambda}\left(\frac{\frac{1}{\sigma_i}+2N^2\lambda^2\sum_{j=1}^N\sigma_ju_j^2-\sigma_iu_i^2N^2\lambda^2}{\sum_{j=1}^N\sigma_ju_j^2}\right)>0.
\end{align*}
With the growth behavior of $\lambda\sigma_i$ it is easy to deduce the limit behaviors.
\qed

Thus $x_i^*$ behaves inversely and monotonically with respect to changes in $\sigma_i$.

In the previous section we saw that for $\sigma_j\equiv \sigma>0$, the \textbf{HMO} is always achieved at the compromise state found on the $N$-sphere determined by $\sigma,u_i$. If we only change the first stubbornness parameter $\sigma_1$, then the \textbf{HMO} state is found on the $N$ dimensional hyperellipsoid defined by
\begin{align*}
    \sigma_1(u_1-x_1)x_1+\sigma\sum_{j=2}^N(u_j-x_j)x_j=0.
\end{align*}
Varying $\sigma_1$ provides insight on when the \textbf{HMO} is found at the compromise state, or when we need to prune and/or strengthen the network to produce polarization or consensus.
\begin{theorem}\label{t:prune1}
    Let $\sigma_1=\mu$ and $\sigma_j=\sigma$ for all $j=2\ldots N$. Then there exists values $\mu_-,\mu_+$ such that $0<\mu_-<\sigma<\mu_+$ where
    \begin{enumerate}
        \item the \textbf{HMO} is achieved at the compromise state via the algorithm produced in Theorem \ref{algorithm} for all $\mu \in (\mu_-,\mu_+)$,
        \item For $\mu \in [\mu_+,\infty)$, the first agent $x_1$ should be pruned from the matrix yielding a new \textbf{HMO(p)} with $x_1^*=u_1$ and $x_j^*$ found on the $(N-1)$-sphere determined by Theorem \ref{t:uniform}.
        
        \item As $\mu$ decreases in $(0,\mu_-)$ we provide an algorithm for pruning the connectivity matrix where successively each $x_2$, then $x_3$, and so on are pruned until only $x_1$ and $x_N$ remain with a final \textbf{HMO(p)} achieved at a consensus between $x_1^*$ and $x_N^*$ and all other agents polarized $x_j^*=u_j$, for $j=2\ldots N-1$.
    \end{enumerate}  
\end{theorem}

\Proof
If $\mu=\sigma$, then we are at the homogeneous setting where the \textbf{HMO} is attained at the compromise state. Recall that $g_j=\sigma_j(x_j^*-u_j)x_j^*$ so that $g_1=\mu(x_1^*-u_1)x_1^*$ and $g_j=\sigma(x_j^*-u_j)x)j^*$ for $j=2\ldots N$. Since we know from Theorem \ref{t:uniform} that at $\mu=\sigma$, the criterion \eqref{suf} holds, by continuity, there exists a $\varepsilon>0$ small enough such that the criterion \eqref{suf} must hold on a small interval $\mu \in (\sigma-\varepsilon,\sigma+\varepsilon)$. However, we can find exactly what these values are.

First, by Lemma \ref{parameter} we see that $x_1^*$ decays as $\mu$ increases. Since $\sigma_j\equiv \sigma_j$ for $j=2 \dots N$ we can also see that
\begin{align*}
    \frac{d x_j^*}{d\mu}=-\frac{1}{2N\lambda^2\sigma}\frac{d\lambda}{d\mu}, \ \ \text{for all} \ j=2 \ldots N.
\end{align*}
Therefore varying $\mu$ not only forces a monotonic behavior on $x_1^*(\mu)$, but implies that the shift of all other $x_j^*(\mu)$ is uniform. Hence, for all $j,k\neq 1$ we have
\begin{align*}
    x_{j}^*-x_k^*&=\frac{1}{2}(u_j-u_k),\\
    x_1^*-x_j^*&=\frac{1}{2}(u_1-u_2)+\frac{1}{2N\lambda}\left(\frac{1}{\mu}-\frac{1}{\sigma}\right).
\end{align*}
Therefore for $\mu>\sigma$, there can be no reordering of agents. Thus for $\mu>\sigma$ the conditions $x_1^*<\ldots <x_N^*$ and $u_1<\ldots <u_N$ must hold.

Now let us define $\mu_k$ for each $k=1\ldots N-1$ such that
\begin{align*}
    \mu_k(x_1^*-u_1)x_1^*+\sigma\sum_{j=2}^k(x_j^*-u_j)x_j^*=0.
\end{align*}
Plugging in for $x_j^*$ and $\lambda$ yields $\mu_k$ satisfying the following quadratic equations
\begin{align}\label{muk}
    f_k(\mu_k):=\mu_k^2\frac{u_1^2}{\sigma}\left(k-N\right)&+\mu_k\left((k-1)\sum_{j=2}^Nu_j^2-(N-1)\sum_{j=2}^ku_j^2\right)\nonumber\\
    & \ \ \ \ \ \ +\sigma\left(\sum_{j=2}^Nu_j^2-\sum_{j=2}^ku_j^2\right)=0, \ \ k=1\dots N-1
\end{align}

Plugging $\mu_k=\sigma$ into \eqref{muk} yields
\begin{align*}
    f_k(\sigma)=k\sum_{j=1}^Nu_j^2-N\sum_{j=1}^ku_j^2>0.
\end{align*}
where the above is positive for each $k$ because of the increasing nature of each $u_j<u_{j+1}$.
As the leading coefficient of $f_k$ is negative for all $k=1\ldots N-1$ we have a unique $\mu_k>\sigma$ satisfying \eqref{muk} for each $k=1 \dots N-1$. Further, we can see that $\mu_1$ which is the value for which $x_1^*=u_1$ is the minimal $\mu_k$. Indeed,
\begin{align*}
    \mu_1=\sigma\sqrt{\frac{\sum_{j=2}^Nu_j^2}{(N-1)u_1^2}},
\end{align*}
while we have the bound
\begin{align*}
    \mu_k\geq \sigma\sqrt{\frac{\sum_{j=k+1}^Nu_j^2}{(N-k)u_1^2}}\geq \mu_1,
\end{align*}
where the lower bound is obtained again from the fact that the linear coefficient in the quadratic is positive due to $u_j<u_j+1$ for all $j$.
 and hence we can define
 $\mu_+:=\min_{k=1\dots N-1}\mu_k=\mu_1>\sigma$. For which the compromise state with adjacency matrix constructed in Theorem \ref{algorithm} is guaranteed to provide the \textbf{HMO} on $[\sigma,\mu_+)$.\\

 Since $\mu_+=\mu_1$, for any $\mu \in [\mu_+,\infty)$ we have $x_1^*\leq u_1$ and hence we prune $x_1$ from the network so that the new \textbf{HMO(p)} is given by
 \begin{align*}
     x_1^*=u_1, \ \ \ x_j^*=\frac{1}{2(N-1)\lambda_p\sigma}+\frac{u_j}{2}, \ \ \ \lambda_p=\frac{1}{\sigma\sqrt{(N-1)\sum_{j=2}^Nu_j^2}}.
 \end{align*}
 With $x_1$ pruned from the network the remainder of $\sigma_j=\sigma$, $j=2\ldots N$ and by Theorem \ref{t:uniform} the \textbf{HMO(p)} is achieved with the prescribed adjacency matrix therein.\\

 To achieve the lower region where $\mu \in (\mu_-,\sigma)$, we have again from Lemma \ref{parameter} that $x_1^*(\mu)$ increases as $\mu$ decreases, hence there is no value for which $x_1^*(\mu)=u_1$ for $\mu<\sigma$. However, there due to the asymptotic behavior of each $x_j^*(\mu)\to \frac{u_j}{2}$ as $\mu\to 0$, there are strictly decreasing finite values $\mu_{x_k}$ such that $x_1^*(\mu_{x_k})=x_k^*(\mu_{x_k})$, for each $k=2,...,N$ which could pose problematic for the construction of the adjacency matrix. Indeed, after each $\mu_{x_k}$, the agents must be re-ordered in order to be amenable to the construction algorithm in Theorem \ref{algorithm}. When this occurs, the analysis provided for the quadratic equations $f_k$ given in \eqref{muk} breaks down since the monotonicity of $u_j$ fails upon re-ordering. Indeed we have
 \begin{align*}
     u_1<\ldots <u_n, \ \ &\mu \in (\mu_{x_2},\infty),\\
     u_2<u_1<u_3<\ldots <u_n, \ \ &\mu \in (\mu_{x_3}, \mu_{x_2}),\\
     u_2<u_3<\ldots<u_k<u_1<u_{k+1}<\ldots < u_n, \ \ &\mu \in (\mu_{x_{k+1}},\mu_{x_k}).
 \end{align*}
We can explicitly compute these values by solving the equations $x_1^*=x_k^*$, leading to the following cubic equations
 \begin{align*}
     g_k(\mu_{x_k})&:=(\mu_{x_k})^3\frac{u_1^2}{\sigma^2}+\frac{(\mu_{x_k})^2}{\sigma}\left(\sum_{j=2}^Nu_j^2-2u_1^2-(N-1)(u_k-u_1)^2\right)+\\
     & \ \ \ \ \ \ \ \ +\mu_{x_k}\left(-2\sum_{j=2}^Nu_j^2+u_1^2-(u_k-u_1)^2\right)+\sigma\sum_{j=2}^Nu_j^2=0
 \end{align*}
 Checking at $\mu^k=0$ and $\mu^k=\sigma$ yields
 \begin{align*}
     g_k(0)=\sigma\sum_{j=2}^Nu_j^2>0, \ \ \  g_k(\sigma)=-\sigma N(u_k-u_1)^2<0,
 \end{align*}
 and hence for each $k=2,...,N$ there is a unique solution to $g(\mu_{x_k})=0$ in the interval $(0,\sigma)$. As $\mu_{x_2}$ is guaranteed to be the largest, we have proved the first part of the theorem at least that on $(\mu_{x_2},\mu_+)$. Note that the \textbf{HMO} is not guaranteed to fail exactly at $\mu_{x_2}$. Indeed, if at $\mu=\mu_{x_2}$ if $x_2^*>u_2$, then after re-ordering there is no violation of the algorithm from Theorem \ref{algorithm}.

 Thus we need also to find the values $\mu_{u_k}$ for which $x_k^*=u_k$. Solving this yields the following quadratic equations
 \begin{align*}
     h_k(\mu_{u_k}):=(\mu_{u_k})^2\left(\frac{u_1^2}{\sigma^2u_k^2}\right)+\mu_{u_k}\left(\frac{1}{\sigma u_k^2}\sum_{j=2}^Nu_j^2-\frac{N-1}{\sigma}\right)-1=0.
 \end{align*}
Where we see by checking at $h_k(\sigma)$, that there is a unique solution $\mu_{u_k}\in (0,\sigma)$ if and only if $\frac{\sum_{j=1}^Nu_j^2}{u_k^2}>N$. Otherwise there is a unique $\mu_{u_k}\in [\sigma,\infty)$. Therefore $\mu_-=\min (\mu_{x_2},\mu_{u_2})<\sigma$.\\

Now for any $\mu \in (0,\mu_-]$ the first step is to prune $x_2$ from the network giving
\begin{align}
    &x_1^*=\frac{1}{2(N-1)\lambda_{p_1}\mu}+\frac{u_1}{2}, \ \ x_2^*=u_2, \ \ x_j^*=\frac{1}{2(N-1)\lambda_p\sigma}+\frac{u_j}{2}, \label{HMOp1-1}\\
    & \ \ \ \ \ \ \ \ \ \ \ \ \ \ \ \lambda_{p_1}=\frac{1}{N-1}\sqrt{\frac{\frac{1}{\mu}+\frac{N-2}{\sigma}}{\mu u_1^2+\sigma\sum_{j=3}^Nu_j^2}}.\label{HMOp1-2}
\end{align}
Starting from $\mu=\mu_-$ we continue by decreasing $\mu$ again. It is important to note that it is possible after pruning that $x_3^*<x_1^*$ and $x_3^*<u_3$ in which case, before varying $\mu<\mu_-$ we must prune $x_3$ from the network. If this is not the case and $u_1<x_1^*<x_3^*<x_N^*<u_N$, then the new \textbf{HMO($p_1$)} is given by \eqref{HMOp1-1}-\eqref{HMOp1-2}. Then we can continue reducing $\mu$ and subsequently pruning $x_3 \ldots x_{N-1}$. Finally we have reduced to the two dimensional case comparing $x_1^*$ and $x_N^*$, the final step for small $\mu$ is then not to prune, but to strengthen the connection between $x_1$ and $x_N$ to get consensus between these two agents, so that \textbf{HMO($p_{\text{final}}$)} is given by
\begin{align*}
    x_1^*=x_N^*=\frac{\frac{\mu}{\sigma}u_1+u_N}{\frac{\mu}{\sigma}+1}, \ \ x_j^*=u_j, \ \ j=2\ldots N-1.
\end{align*}
\qed

We chose the stubbornness parameter $\sigma_1$ corresponding to $u_1=\min_j u_j$ to vary, however, similar algorithms can be produced for varying any particular $\sigma_j$, so long as the original vectors $(\sigma, u)$ produced an \textbf{HMO} that satisfied condition \eqref{suf}.\\

In the case of given parameters $(\sigma, u)$ such that \eqref{suf}  does not hold, the \textbf{HMO} calculated in \eqref{equ} is not representative of a fixed point for any symmetric nonnegative adjacency matrix $A$. This means that the computed \textbf{HMO} requires that some connections within $A$ be negative, $a_{ij}<0$ for some $i,j$. Therefore, in order to achieve a next best \textbf{HMO}, we must adjust the interactions between agents. In general, after computing $x_j^*$ and re-ordering, the new order of $u_j$ and the relationships between $x_1$ and $u_1$ as well as $x_N$ and $u_N$ indicate what graph modifications should be made.

For instance, if after the computation of $x_j^*$, there is no re-ordering, and \eqref{suf} fails, then one should either prune the first node $k$ for which \eqref{suf} fails (bottom pruning), or doing the algorithm symmetrically from the top the first $k$ for which the negative of \eqref{suf} fails (top pruning). This pruning cannot affect the ordering and thus one can continue choosing between bottom and top pruning appropriately until a suitable smaller network structure remains for which \eqref{suf} holds.

On the other hand, if for example a complete re-ordering occurs so that $u_1>\ldots >u_N$, then rather than pruning the new \textbf{HMO} occurs at the consensus value achieved by uniformly strengthening all network connections towards infinity.

\section{Conclusions and remarks}\label{sect8}

We have analyzed a recently proposed nonlinear model of opinion dynamics that incorporates both peer influence and individual conviction across a networked community. Our central focus was the identification and characterization of the highest mean opinion (HMO) state, and its relationship to consensus in heterogeneous populations.

Our findings include:

\begin{itemize}
\item Characterization of consensus: in our model, perfect consensus is only achievable when all agents share identical conviction values. In realistic heterogeneous settings, variance in equilibrium opinions persists, but the equilibrium approaches consensus as connectivity increases.

\item HMO does not imply consensus: through both theoretical analysis and concrete examples, we demonstrated that maximizing the average opinion often leads to heterogeneous states rather than uniform agreement.

\item Optimization of mean opinion: we provided necessary and sufficient conditions to identify and construct the optimal connectivity matrix that yields the {\bf HMO} at a comprise state.

\item Pruning and Strengthening the connectivity matrix: When the conditions to reach the \textbf{HMO} fail due to heterogeneity of stubbornness parameters, we provide guidelines for how to reach a next best \textbf{HMO(p)} via pruning `problematic' agents from the connectivity matrix or strengthening the connection between particular agents.

\item In large communities with subgroup structures, optimal configurations typically involve selective inter-group connections rather than uniform intra-group links. This suggests that tailored connections between different clusters can be more effective than reinforcing internal cohesion.
\end{itemize}

Regarding possible extensions or variations of our study, future work could address the limitations of our current model by considering effects like dynamic conviction values, time-varying networks, or stochastic perturbations. In particular, the first aspect is very natural, because personal conviction values are not always fixed, and in many cases may also be influenced on the opinions of other individuals. To this regard, in the next future we plan to consider a simple extension of our model \eqref{model} by considering the personal conviction values $u_i$ not as fixed parameter but temporal functions, according to a rule
$$
u_i'(t)=f_i(x,u_i).
$$
For instance, a natural law that can be proposed is 
\begin{equation}\label{convict}
u_i'(t)=a_i\left({\cal M}(x)-u_i\right),
\end{equation}
with $a_i\geq 0$ for $i=1,\ldots,N$. For $a_i>0$, it means that the personal conviction value of each agent at a certain moment is attracted to the mean opinion value. Our conjecture is that system \eqref{model}-\eqref{convict} always will tend to a consensus state, but to understand what would be the optimal connectivity matrix $A$ to reach the {\bf HMO} value seems a challenging problem. A correct answer to this and many other related questions would contribute to further enriching our understanding of collective opinion formation.

\end{document}